\documentclass[smallextended]{svjour3}       
\smartqed  
\usepackage{graphicx}
\usepackage{subfigure}
\usepackage{amssymb,amsmath}
\usepackage{natbib}

\journalname{}
\begin{document}

\title{Optimal control as a graphical model inference problem}



\author{Hilbert J. Kappen \and Vicen\c{c} G\'omez \and Manfred Opper}

\institute{Hilbert J. Kappen \at
        Donders Institute for Brain Cognition and Behaviour\\
        Radboud University Nijmegen\\
        6525 EZ Nijmegen, The Netherlands\\
              \email{b.kappen@science.ru.nl}           
           \and
        Vicen\c{c} G\'omez \at
        Donders Institute for Brain Cognition and Behaviour\\
        Radboud University Nijmegen\\
        6525 EZ Nijmegen, The Netherlands\\
              \email{v.gomez@science.ru.nl}           
           \and
        Manfred Opper \at
        Department of Computer Science\\
        D-10587 Berlin, TU Berlin, Germany\\
               \email opperm@cs.tu-berlin.de
}

\date{}

%
%
\newcommand{\HRule}{\rule{\linewidth}{2mm}}

\newcommand{\bmpt}[1]{\begin{minipage}[t]{#1\textwidth}}
\newcommand{\bmp}[1]{\begin{minipage}{#1\textwidth}}
\newcommand{\nck}[2]{\left(\begin{tabular}{c}#1\\#2\end{tabular}\right)}
\newcommand{\emp}{\end{minipage}}
\newcommand{\Na}{\mathsf{Na}^+}
\newcommand{\K}{\mathsf{K}^+}
\newcommand{\Cl}{\mathsf{Cl}^-}
\newcommand{\Ca}{\mathsf{Ca}^{2+}}
\newcommand{\no}{\nonumber}
\newcommand{\n}{\noindent}
\newcommand{\On}{{\cal O}(n)}
\newcommand{\Osqrtn}{{\cal O}(\sqrt{n})}
\newcommand{\Ai}{\mathrm{Ai}}
\newcommand{\Bi}{\mathrm{Bi}}
\newcommand{\bs}{{\bf s}}
\newcommand{\cR}{{\cal R}}
\newcommand{\cM}{{\cal M}}
\newcommand{\cP}{{\cal P}}
\newcommand{\cS}{{\cal S}}
\newcommand{\cC}{{\cal C}}
\newcommand{\cQ}{{\cal Q}}
\newcommand{\cT}{{\cal T}}
\newcommand{\cN}{{\cal N}}
\newcommand{\cA}{{\cal A}}
\newcommand{\cV}{{\cal V}}
\newcommand{\cO}{{\cal O}}
\newcommand{\bv}{\begin{verbatim}}
\newcommand{\ev}{\end{verbatim}}
\newcommand{\Ei}{E_{\mbox{\small{int}}}}
\newcommand{\ben}{\begin{enumerate}}
\newcommand{\een}{\end{enumerate}}
\newcommand{\be}{\begin{equation}}
\newcommand{\ee}{\end{equation}}
\newcommand{\bea}{\begin{eqnarray}}
\newcommand{\eea}{\end{eqnarray}}
\newcommand{\beaa}{\begin{eqnarray*}}
\newcommand{\eeaa}{\end{eqnarray*}}
\newcommand {\bd}{\begin{description}}
\newcommand {\ed}{\end{description}}
\newcommand{\bi}{\begin{itemize}}
\newcommand{\ei}{\end{itemize}}
\newcommand{\bc}{\begin{center}}
\newcommand{\ec}{\end{center}}
\newcommand {\bq}{\begin{quotation}}
\newcommand {\eq}{\end{quotation}} 
\newcommand {\bdo}{\begin{document}} 
\newcommand {\edo}{\end{document}} 
\newcommand{\eql}{\: = \:}
\newcommand{\dfl}{\: \equiv \:}
\newcommand{\onder}[2]{#1_{\mbox{\scriptsize #2}}}
\newcommand {\boven}[2]{#1^{\mbox{\scriptsize #2}}}
\newcommand{\argmin}[1]{\mathop{\rm argmin}\limits_{#1}}
\newcommand{\klein}[2]{#1_{\mbox{\tiny #2}}}
\newcommand{\doo}{\mbox{do}}
\newcommand{\fracd}[2]{{{\displaystyle #1} \over {\displaystyle #2}}}

%
%
\newcommand{\intarg}[2]{\!#2\,d#1}
\newcommand{\del}[3]{\left#1 #3 \right#2}
\newcommand{\av}[1]{\del{<}{>}{#1}}
\newcommand{\bra}[1]{\del{<}{|}{#1}}
\newcommand{\ket}[1]{\del{|}{>}{#1}}
\newcommand{\avmini}[1]{\del{<}{>_{\setminus i}}{#1}}
\newcommand{\avminij}[1]{\del{<}{>_{\setminus i,j}}{#1}}
\newcommand{\sminij}{s_{\setminus i,j}}
\newcommand{\ui}{\avmini{u_i^2}}
\newcommand{\hi}{\avmini{h_i}}
\newcommand{\br}[1]{\del{\{}{\}}{#1}}
\newcommand{\ar}[1]{\del{(}{)}{#1}}
\newcommand{\bl}[1]{\del{[}{]}{#1}}

\newcommand {\Hxx}{H^{(xx)}}
\newcommand {\Hxy}{H^{(xy)}}
\newcommand {\Hyx}{H^{(yx)}}
\newcommand {\Hyy}{H^{(yy)}}
\newcommand {\Cxx}{C^{(xx)}}
\newcommand {\Cxy}{C^{(xy)}}
\newcommand {\Cyx}{C^{(yx)}}
\newcommand {\Cyy}{C^{(yy)}}
\newcommand {\cx}{c^{(x)}}
\newcommand {\hy}{h^{(y)}}
\newcommand{\bt}{\bar{t}}
\newcommand {\bx}{\bar{x}}
\newcommand {\bX}{\bar{X}}
\newcommand {\by}{\bar{y}}
\newcommand {\bw}{\bar{w}}
\newcommand {\bh}{\bar{h}}
\newcommand {\vs}{\vec{s}}
\newcommand {\vbfs}{{\bf s}}
\newcommand {\vv}{\vec{v}}
\newcommand {\vh}{\vec{h}}
\newcommand {\vk}{\vec{k}}
\newcommand {\vz}{\vec{z}}
\newcommand {\vx}{\vec{x}}
\newcommand {\valpha}{\vec{\alpha}}
\newcommand {\vm}{\vec{m}}
\newcommand {\vtheta}{\vec{\theta}}
\newcommand {\vgamma}{\vec{\gamma}}
\newcommand {\qd}{q_{\partial_i}}
\newcommand {\qm}{q_{\setminus i}}
\newcommand{\vxv}{\vec{x}_v}
\newcommand{\vxh}{\vec{x}_h}
\newcommand {\vy}{\vec{y}}
\newcommand {\bfy}{{\bf y}}
\newcommand {\bfx}{{\bf x}}
\newcommand {\bft}{{\bf t}}
\newcommand {\bfm}{{\bf m}}
\newcommand {\bfu}{{\bf u}}
\newcommand {\bfw}{{\bf w}}
\newcommand {\bfW}{{\bf W}}
\newcommand {\bfmu}{{\bf \mu}}
\newcommand {\vw}{\vec{w}}
\newcommand {\vb}{\vec{b}}
\newcommand {\vu}{\vec{u}}
\newcommand {\vW}{\vec{W}}
\newcommand {\vM}{\vec{M}}
\newcommand {\sign}{\mbox{sign}}
\newcommand {\sech}{\mbox{sech}}
\newcommand {\hmu}{h^{\mu}_j}
\newcommand {\gmu}{g^{\mu}_j}
\newcommand {\gjkmu}{g^{\mu}_{jk}}
\newcommand {\sia}{s_i^\alpha}
\newcommand {\sjb}{s_j^\beta}
\newcommand {\hia}{h_i^\alpha}
\newcommand {\hjb}{h_j^\beta}
\newcommand {\mia}{m_i^\alpha}
\newcommand {\mib}{m_i^\beta}
\newcommand {\mjb}{m_j^\beta}
\newcommand {\wijab}{w_{ij}^{\alpha\beta}}
\newcommand {\chiijab}{\chi_{ij}^{\alpha\beta}}
\newcommand {\Cijab}{C_{ij}^{\alpha\beta}}
\newcommand {\Hijab}{H_{ij}^{\alpha\beta}}
\newcommand {\he}{h_{\mbox{ext}}}
\newcommand {\xia}{x_{i,\alpha}}
\newcommand {\vxk}{\vec{x}_{\mbox{{\small known}}}}
\newcommand {\vxu}{\vec{x}_{\mbox{{\small unknown}}}}
\newcommand {\Pab}{P_{\alpha\beta}}
\newcommand {\dab}{\delta^{\alpha\beta}}
\newcommand {\dij}{\delta_{ij}}
\newcommand {\erf}{\mbox{Erf}}

\maketitle

\begin{abstract}
We reformulate a class of non-linear stochastic optimal control problems
introduced by \cite{NIPS2006_691} as a Kullback-Leibler (KL) minimization problem.  As a result,
the optimal control computation reduces to an inference computation and
approximate inference methods can be applied to efficiently compute approximate
optimal controls. We show how this KL control theory contains the path integral
control method as a special case.  We provide an example of a block stacking
task and a multi-agent cooperative game where we demonstrate how approximate
inference can be successfully applied to instances that are too complex for
exact computation. We discuss the relation of the KL control approach to other
inference approaches to control. 
\end{abstract}

\section{Introduction}
Stochastic optimal control theory deals with the problem to compute an
optimal set of actions to attain some future goal. With each action
and each state a cost is associated and the aim is to minimize the
total future cost. Examples are found 
in many contexts such as motor control tasks for robotics, planning
and scheduling tasks or managing a financial portfolio. The
computation of the optimal control is typically very difficult due to
the size of the state space and the stochastic nature of the problem.

The most common approach to compute the optimal control is through
the Bellman equation. For the finite horizon discrete time case, this
equation results from a dynamic programming argument that expresses
the optimal cost-to-go (or value function) at time $t$ in terms of the
optimal cost-to-go at time $t+1$. For the infinite horizon case, the
value function is independent of time and the Bellman equation becomes
a recursive equation. In continuous time, the Bellman equation becomes
a partial differential equation.

For high dimensional systems or for continuous systems the state
space is huge and the above procedure cannot be directly applied. A common
approach to make the computation tractable is a function approximation
approach where the value function is parameterized in terms of a number
of parameters \citep{bertsekas-tsitsiklis1996}. Another promising
approach is to exploit graphical structure that is present in the
problem to make the computation more efficient
\citep{boutilier1995,koller1999}. However, this graphical structure
is in general not inherited by the value function, and thus the
graphical representation of the value function may not be appropriate.

In this paper, we introduce a class of stochastic optimal control
problems where the control is expressed as a probability
distribution $p$ over future trajectories given the current state and where the control
cost can be written as a Kullback-Leibler (KL) divergence between $p$ and some interaction
terms.  The optimal control is given by minimizing the KL divergence,
which is equivalent to solving a probabilistic inference problem
in a dynamic Bayesian network.  The optimal control is given in
terms of (marginals of) a probability distribution over future
trajectories.  The formulation of the control problem as an
inference problem directly suggests exact inference methods
such as the Junction Tree method (JT) \citep{lauritzen88}
or a number of well-known approximation
methods, such as the variational method \citep{jordan96}, belief
propagation (BP) \citep{murphy99}, the cluster variation method (CVM) or generalized
belief propagation (GBP) \citep{yedidia00} or Markov Chain Monte Carlo (MCMC)
sampling methods.  We refer to this class of problems as KL control
problems.

The class of control problems considered in this paper is identical
as in \cite{todorov_ieee2008,NIPS2006_691,todorov_pnas2009}, who shows that the
Bellman equation can be written as a KL divergence of probability
distributions between two adjacent time slices and that the Bellman
equation computes backward messages in a chain as if it were an
inference problem.  The novel contribution of the present paper is
to identify the control cost with a KL divergence instead of making this identification in the
Bellman equation. The immediate consequence is that the optimal
control problem {\em is identical to} a graphical model inference problem 
that can be approximated using standard methods. 

We also show how KL control reduces to the previously proposed path
integral control problem \citep{kappen_prl05} when noise is Gaussian
in the limit of continuous space and time. This class of control
problem has been applied to multi-agent problems using a graphical
model formulation and junction tree inference in
\citet{wimw_uai06,wimw_aamas2007} and approximate inference in
\cite{broek2007alamas,broek2006}.  In robotics,
\citet{theodorou_ADPRLAIS_2009,theodorou_aistats2010,theodorou_icra2010} has
shown the the path integral method has great potential for application.
They have compared the path integral method with some
state-of-the-art reinforcement learning methods, showing very
significant improvements. In addition, they have successful implemented
the path integral control method to a walking robot dog. The path
integral approach has recently been applied to the control of
character animation \citep{dasilva2009}.

\section{Control as KL minimization}
\label{method}
Let $x=1,\ldots, N$ be a finite set of states, $x^t$ denotes the state
at time $t$. Denote by $p^t(x^{t+1}|x^t,u^t)$ the
Markov transition probability at time $t$ under control $u^t$ from state
$x^t$ to
state $x^{t+1}$. Let $p(x^{1:T}|x^0,u^{0:T-1})$ denote the probability to
observe the trajectory $x^{1:T}$ given initial state $x^0$ and control
trajectory $u^{0:T-1}$. 

If the system at time $t$ is in state $x$ and takes action $u$ to state
$x'$, there is an
associated cost $\hat{R}(x,u,x',t)$. The control problem is to find the sequence
$u^{0:T-1}$ that minimizes the expected future cost
\begin{align}
C(x^0,u^{0:T-1})&=\sum_{x^{1:T}} p(x^{1:T}|x^0,u^{0:T-1})
\sum_{t=0}^{T} \hat{R}(x^{t},u^{t},x^{t+1},t)\notag\\
&=\left\langle\sum_{t=0}^{T}
\hat{R}(x^{t},u^{t},x^{t+1},t)\right\rangle
\label{Cmdp}
\end{align}
with the convention that $\hat{R}(x^{T},u^{T},x^{T+1},T)=R(x^T,T)$ is the cost
of the final state and $\langle\rangle$ denotes expectation with respect to $p$.
Note, that $C$ depends on $u$ in two ways: through $\hat{R}$ and through the
probability distribution of the controlled trajectories
$p(x^{1:T}|x^0,u^{0:T-1})$.

The optimal control is normally computed using the Bellman equation, which
results from a dynamic programming argument \citep{bertsekas-tsitsiklis1996}. 
Instead, we will consider the restricted class of control problems for which
$C$ in Equation \eqref{Cmdp} can be written as a KL divergence.  As a
particular case, we consider that $\hat{R}$ is the sum of a control dependent
term and a state dependent term.  We further assume the existence of a 'free'
(uncontrolled) dynamics $q^t(x^{t+1}|x^t)$, which can be any first order Markov
process that assigns zero probability to physically impossible state
transitions.

We quantify the control cost as the amount of deviation between
$p^t(x^{t+1}|x^t,u^t)$ and $q^t(x^{t+1}|x^t)$ in KL sense.  Thus,
\begin{align}
\hat{R}(x^t,u^t,x^{t+1},t)=\log
\frac{p^t(x^{t+1}|x^t,u^t)}{q^t(x^{t+1}|x^{t})}+ R(x^t,t)\qquad
t=0,\ldots,T-1
\label{Rspecial}
\end{align}
with $R(x,t)$ an arbitrary state dependent control cost. Equation \eqref{Cmdp}
becomes
\begin{align}
C(x^0,p)&=KL(p||\psi)\notag\\
&=\sum_{x^{1:T}}p(x^{1:T}|x^0)\log\frac{p(x^{1:T}|x^0)}{\psi(x^{1:T}|x^0)}\notag\\
&=KL(p||q)+\langle{R}\rangle\label{C}\\
\psi(x^{1:T}|x^0)&=q(x^{1:T}|x^0) \exp\left(-\sum_{t=0}^{T}
R(x^t,t)\right)\label{psi}
\end{align}

Note, that $C$ depends on the control $u$ only through $p$. Thus, minimizing $C$ 
with respect to $u$ yields: $0=\frac{dC}{du}=\frac{dC}{dp}\frac{dp}{du}$, where the minimization
with respect to $p$ is subject to the normalization constraint
$\sum_{x^{1:T}}p(x^{1:T}|x^0)=1$.  Therefore, a sufficient condition for the optimal control
is to set $\frac{dC}{dp}=0$.
The result of this KL minimization is well known and yields the ``Boltzmann distribution"

\begin{align}
p(x^{1:T}|x^0)&=\frac{1}{Z(x_0)}\psi(x^{1:T}|x^0)
\label{p-optimal}
\end{align}
and the optimal cost
\begin{align}
C(x^0,p)&=-\log
Z(x^0)=-\log\sum_{x^{1:T}}q(x^{1:T}|x^0)\exp\left(-\sum_{t=0}^T
R(x^t,t)\right)
\label{path-integral}
\end{align}
where $Z(x^0)$ is a normalization constant (see Appendix~\ref{appendix:bm}). In
other words, the optimal control solution is the (normalized) product of the
free dynamics and the exponentiated costs. It is a distribution that avoids
states of high $R$, at the same time deviating from $q$ as little as possible.
Note that since $q$ is a first order Markov process, $p$ in
Equation~\eqref{p-optimal} is a first order Markov process as well.

The optimal control in the current state $x^0$ at the current time $t=0$ is given by the marginal probability
\begin{align}
p(x^{1}|x^0)&=\sum_{x^{2:T}}p(x^{1:T}|x^0)\label{marginal}
\end{align}
This is a standard graphical model inference problem, with $p$ given by
Equation~\eqref{p-optimal}.  Since $\psi$ is a chain, we can compute
$p(x^{1}|x^0)$ by backward message passing:
\begin{align*}
\beta^T(x^T) & = 1\\
\beta^t(x^t) & = \sum_{x^{t+1}}\psi_t(x^t,x^{t+1})\beta^{t+1}(x^{t+1})\\
p(x^{t+1}|x^t) & \propto  \psi^t(x^t,x^{t+1})\beta^{t+1}(x^{t+1}).
\end{align*}

\begin{figure}
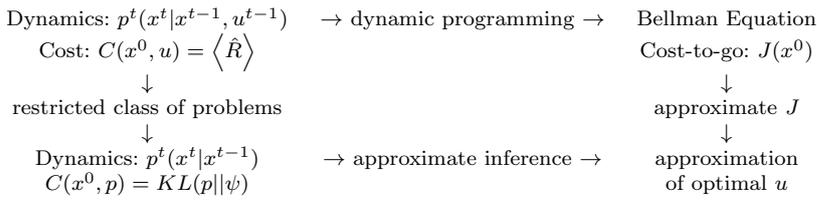

\begin{center}
\begin{tabular}{ccc}
Dynamics: $p^t(x^t|x^{t-1},u^{t-1})$ &$\rightarrow$ dynamic programming $\rightarrow$
&Bellman Equation  \\
Cost: $C(x^0,u)=\left\langle\hat{R}\right\rangle$&&  Cost-to-go: $J(x^0)$\\
$\downarrow$&               &$\downarrow$\\
restricted class of problems&   &approximate $J$\\
$\downarrow$&               &$\downarrow$\\
Dynamics: $p^t(x^t|x^{t-1})$&$\rightarrow$ approximate inference
$\rightarrow$  &approximation\\
$C(x^0,p)=KL(p||\psi)$&&of optimal $u$ 
\end{tabular}
\end{center}
\caption{{\small Overview of the approaches to computing the optimal control.
Top left) The general optimal control problem is formulated as a state transition model
$p$ that depends on the control (or policy) $u$ and a cost $C(u)$
that is the expected $\hat{R}$ with respect to the controlled
dynamics $p$. 
The optimal control
is given by the $u$ that minimizes a cost $C(u)$.
Top right) The traditional approach is to introduce
the notion of cost-to-go or value function $J$, which satisfies the
Bellman equation. The Bellman equation is derived using a dynamic
programming argument.
Bottom right) For large problems, an approximate representation of $J$ is
used to solve the Bellman equation which yields the optimal control.
Bottom left) The
approach in this paper is to consider a class of control problems for
which $C$ is written as a KL divergence. The computation
of the optimal control (optimal $p$) becomes a statistical inference problem, that can be approximated using
standard approximate inference methods.}}
\label{approach}
\end{figure}

The interpretation of the Bellman
equation as message passing for the KL control problems was first established
in \cite{todorov_ieee2008}.
The difference
between the KL control computation and the standard computation
using the Bellman equation is schematically illustrated in
Figure \ref{approach}.

The optimal cost, Equation \eqref{path-integral}, is minus the log partition
sum and is the expectation value of the exponentiated state costs $\sum_{t=0}^T
R(x^t,t)$ under the {\em uncontrolled} dynamics $q$. This is a surprising
result, because it means that we have a closed form solution for the optimal
cost-to-go $C(x^0,p)$ in terms of the known quantities $q$ and $R$.

A result of this type was 
previously obtained in \cite{kappen_prl05} for a class of continuous non-linear
stochastic control problems. Here, we show that a slight
generalization of
this problem ($g_{ai}(x,t)=1$ in \cite{kappen_prl05}) is obtained as a special case of the 
present KL control formulation.
Let $x$ denote an $n$-dimensional real vector with components $x_i$. We
define the stochastic dynamics 
\begin{align}
dx_i&=f_i(x,t)dt + \sum_{a} g_{ia}(x,t)(u_a dt+ d\xi_a) \label{dynamics_pi}
\end{align}
with $f_i$ an arbitrary function, $d\xi_a$ an $m$-dimensional Gaussian process with covariance
matrix $\left\langle d\xi_a d\xi_b\right\rangle=\nu_{ab} dt$ and
$u_a$ an $m$-dimensional
control vector. 
The distribution over trajectories 
is given by
\begin{align}
p(x^{dt:T}|x^0,
u^{0:T-dt})&=\prod_{s=0}^{T-dt}\mathcal{N}(x^{s+dt}|x^s+(f^s +g^s u^s) dt ,g^s \nu (g^s)^T dt)
\label{controlleddynamics}
\end{align}
with $f^t=f(x^t,t)$ and the distribution over
trajectories under the uncontrolled dynamics is defined as
$q(x^{dt:T}|x^0)=p(x^{dt:T}|x^0, u^{0:T-dt}=0)$.

For this particular choice of $p$ and $q$, the control cost in Equation
\eqref{C} becomes (see Appendix~\ref{appendix:continuous} for a derivation)
\begin{align}
C(x,u(t\rightarrow T))&=\left\langle\phi(x(T))+\int_t^T ds
\frac{1}{2}u(x(s),s)^T\nu^{-1} u(x(s),s)+ R(x(s),s)\right\rangle
\label{cost_pi}
\end{align}
where $\langle\rangle$ denotes expectation with respect to the controlled
dynamics $p$, where the sums become integrals and where we have defined
$\phi(x)=R(x,T)$.  

Equations \eqref{dynamics_pi} and \eqref{cost_pi} define a stochastic
optimal control problem. 
The solution for the optimal cost-to-go for this
class of control problems can be shown to be given as a so-called path
integral, an integral over trajectories, which is the continuous time
equivalent of the sum over trajectories in Equation \eqref{path-integral}.
Note, that the cost of control is quadratic in $u$, but of a
particular form with the matrix $\nu^{-1}$ in agreement with
\cite{kappen_prl05}. 
Thus, the KL
control theory contains the path integral control method as a
particular limit.
As is shown in \cite{kappen_prl05}, this class of
problems admits a solution of the optimal cost-to-go as an integral
over paths, which is similar to Equation \eqref{path-integral}.

\subsection{Graphical model inference}
\label{graphical}
In typical control problems, $x$ has a modular structure 
with components $x=x_1,\ldots,x_n$.  For instance, for a multi-joint
arm, $x_i$ may denote the state of each joint. For a multi-agent
system, $x_i$ may denote the state of each agent. In all such
examples, $x_i$ itself may be a multi-dimensional state vector.
In such cases, the optimal control computation, Equation \eqref{marginal},
is intractable.
However, the following assumptions are likely to be true:
\begin{itemize}
\item The uncontrolled dynamics factorizes over components 
$$q^t(x^{t+1}|x^t)=\prod_{i=1}^n q_i^t(x^{t+1}_i| x^t_i)$$
\item The
interaction between components has a (sparse) graphical structure
$R(x,t)=\sum_{\alpha} R_\alpha(x_\alpha,t)$
with $\alpha$ a subset of the indices $1,\ldots,n$ and $x_\alpha$ the
corresponding variables.
\end{itemize}
Typical examples are multi-agent systems and robot arms. In both cases the
dynamics of the individual components (the individual agents and the different
joints, respectively) are independent {\em a priori}. It is only through the
execution of the task that the dynamics become coupled.

Thus, $\psi$ in Equation \eqref{psi} has a graphical structure that we can
exploit when computing the marginals in Equation \eqref{marginal}. For instance,
one may use the junction tree (JT) method, which can be more efficient than
simply using the backward messages.  Alternatively, we can use any of a large
number of approximate graphical model inference methods to compute the optimal
control.  In the following sections, we will illustrate this idea by applying
several approximate inference algorithms in two different tasks.
\section{Stacking blocks (KL-blocks-world)}
\label{sec:blocksworld}
\begin{figure}
\bc
\includegraphics[width=0.4\textwidth]{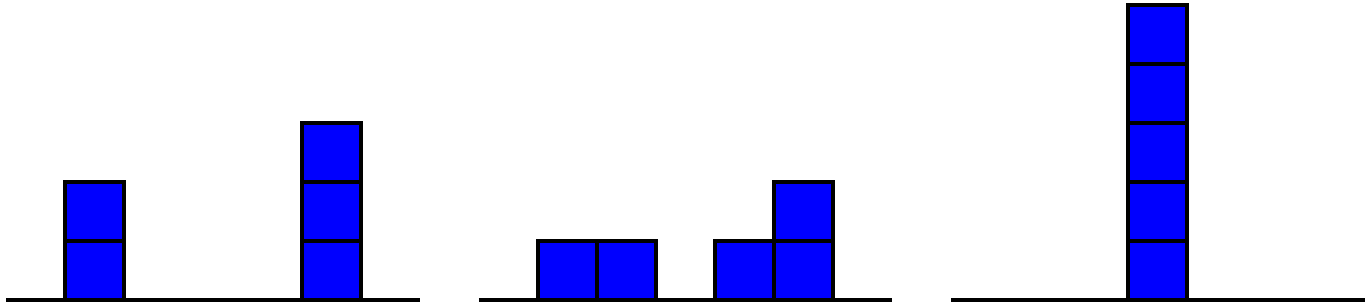}
\caption{Block stacking problem: the objective can be (but is not restricted
to) to stack the initial block configuration (left) into a single stack (right)
through a sequence of single block moves to adjacent positions (middle).}
\label{blocksworld}
\ec
\end{figure}

Consider the example of piling blocks into a tower. This is a classic AI
planning task \citep{russell_norvig03}.  It will be instructive to see how a
variant of this problem is solved as a stochastic control problem, As we will
see, the optimal control solution will in general be a mixture over several
actions.  We define the KL-blocks-world problem in the following way: let there
be $n$ possible block locations on the one dimensional ring (line with periodic
boundaries) as in Figure~\ref{blocksworld},
%
and let $x^t_i\ge 0, i=1,\ldots,n, t=0,\ldots,T$ denote the height of
stack $i$
at time $t$. Let $m$ be the total number of blocks.

At iteration $t$, we allow to move one block from location $k^t$ and
move it to
a neighboring location $k^t+l^t$ with $l^t=-1,0,1$ (periodic boundary
conditions).  Given $k^t,l^t$ and the old state $x^{t-1}$, the new
state is
given as 
\begin{align}
x^{t}_{k^t}&=x^{t-1}_{k^t}-1\\
x^{t}_{k^t+l^t}&=x^{t-1}_{k^t+l^t}+1
\end{align}
and all other stacks unaltered.  We use the uncontrolled distribution
$q$ to
implement these allowed moves. For the purpose of memory efficiency,
we
introduce auxiliary variables $s^t_i=-1,0,1$ that indicate whether
the stack
height $x_i$ is decremented, unchanged or incremented, respectively.
The
uncontrolled dynamics $q$ becomes
$q(k^t)= {\cal U}(1,\ldots,n)$, $q(l^t)= {\cal U}(-1,0,+1)$, 
\begin{align*}
q(s^t|k^t,l^t)  &= \prod_{i=1}^n q(s^t_i|k^t,l^t)\\
q(s^t_i|k^t,l^t)&=
\begin{cases}
\delta_{s^t_i,-1} & \text{for $k^t=i,l^t=\pm 1$}\\
\delta_{s^t_i,+1} & \text{for $k^t+l^t=i,l^t=\pm 1$}\\
\delta_{s^t_i,0}  & \text{otherwise}
\end{cases},
\end{align*}
where $\cal U(\cdot)$ denotes the uniform distribution.
The transition from $x^{t-1}$ to $x^t$ is a
mixture over the values of $k^t,l^t$:
\begin{align}
q(x^{t}|x^{t-1})&=\sum_{k^t,l^t}\prod_{i=1}^n
q(x^{t}_i|x^{t-1}_i,k^t,l^t)q(k^t)q(l^t)\label{blockworld1}\\
q(x^{t}_i|x^{t-1}_i,k^t,l^t)&=\sum_{s_i^t}
q(x^{t}_i|x^{t-1}_i,s^t_i)q(s^t_i|k^t,l^t)\\
q(x^{t}_i|x^{t-1}_i,s^t_i)&= \delta_{x_i^{t},x_i^{t-1}+s^t_i}
\end{align}
Note, that there are combinations of $x_i^{t-1}$ and $s_i^t$ that are
forbidden: we cannot remove a block from a stack of size zero ($x_i^{t-1}=0$
and $s_i^t=-1$) and we cannot move a block to a stack of size $m$
($x_i^{t-1}=m$ and $s_i^t=1$). If we restrict the values of $x^{t}_i$ and
$x^{t-1}_i$ in the last line above to $0,\ldots,m$ these combinations are
automatically forbidden.
\begin{figure}
\bc
\includegraphics[width=0.7\textwidth]{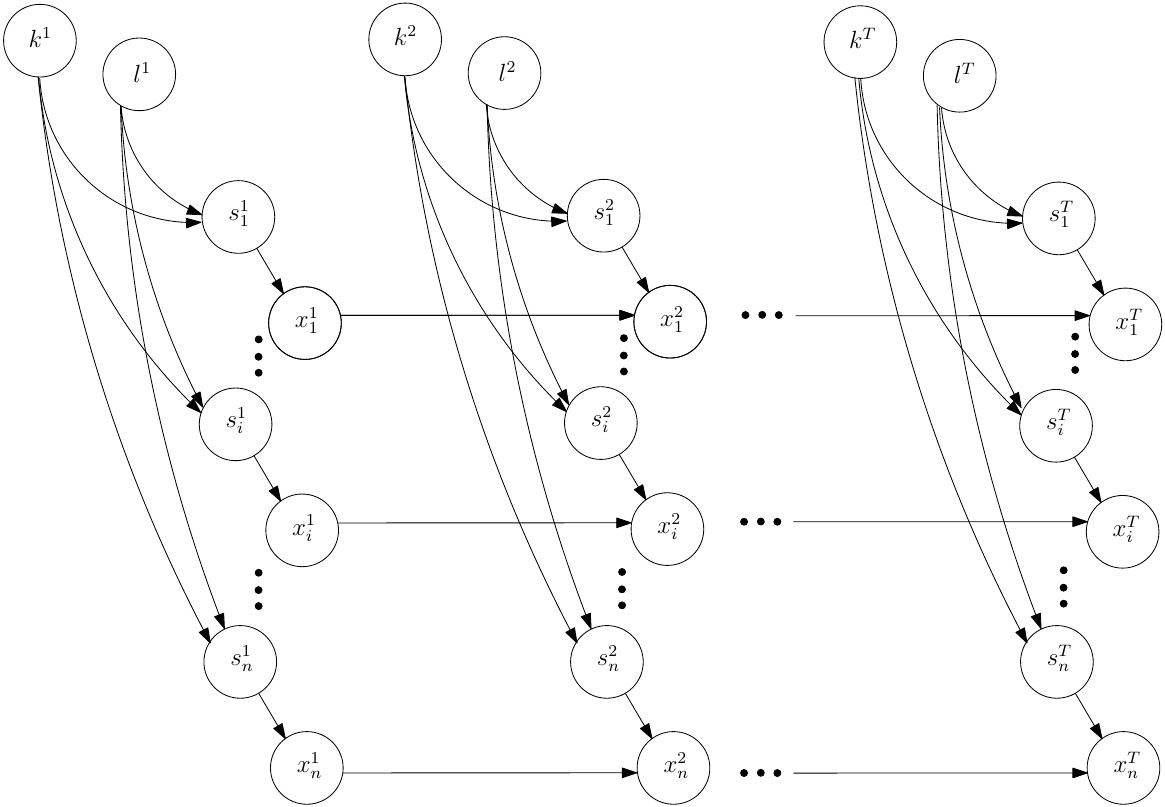}
\caption{
Block stacking problem: Graphical model representation as a dynamic Bayesian
network.  Time runs horizontal and stack positions vertical. At each time, the
transition probability of $x^t$ to $x^{t+1}$ is a mixture over the variables
$k^t,l^t$.  The initial state is ``clamped" to a given configuration by
conditioning on the variables $x^1$.  To force a goal state or final
configuration, the final state $x^T$ can also be ``clamped" (see Section
\ref{sec:KLblocks-goal}).
}
\label{graph}
\ec
\end{figure}

Figure~\ref{graph} shows the graphical model
associated with this representation.  Notice that the graphical structure for
$q$ is efficient compared to the naive implementation of $q(x^t|x^{t-1})$ as a
full table. Whereas the joint table requires $m^{n}$ entries, the graphical
model implementation requires $Tn$ tables of sizes $n\times 3\times 3$ 
for $p(s^t|k^t,l^t)$ and $n\times n\times 3$ for $p(x^t|x^{t-1},s^t)$.
%
In addition, the graphical structure can be
exploited by efficient approximate inference methods.

Finally, a possible state cost can be defined as the entropy of the
distribution of blocks: 
\begin{align}
R(x)&=-\lambda \sum_i \frac{x_i}{m} \log \frac{x_i}{m},
\label{eq:R}
\end{align}
with $\lambda$ a positive number to indicate the strength.  Since
$\sum_i x_i$
is constant (no blocks are lost), the minimum entropy solution puts all blocks
on one stack (if enough time is available).  The control problem is to find the
distribution $p$ that minimizes $C$ in Equation \eqref{C}. 
\subsection{Numerical results}
In the next section, we consider two particular problems.  First, we are
interested in finding a sequence of actions that, starting in a given initial
state $x^0$, reach a given goal state $x^T$, without state cost.  Then we
consider the case of entropy minimization, with no defined goal state and
nonzero state cost.
\subsubsection{Goal state and $\lambda = 0$}
\label{sec:KLblocks-goal}
\begin{figure}
\bc
\includegraphics[width=0.5\textwidth]{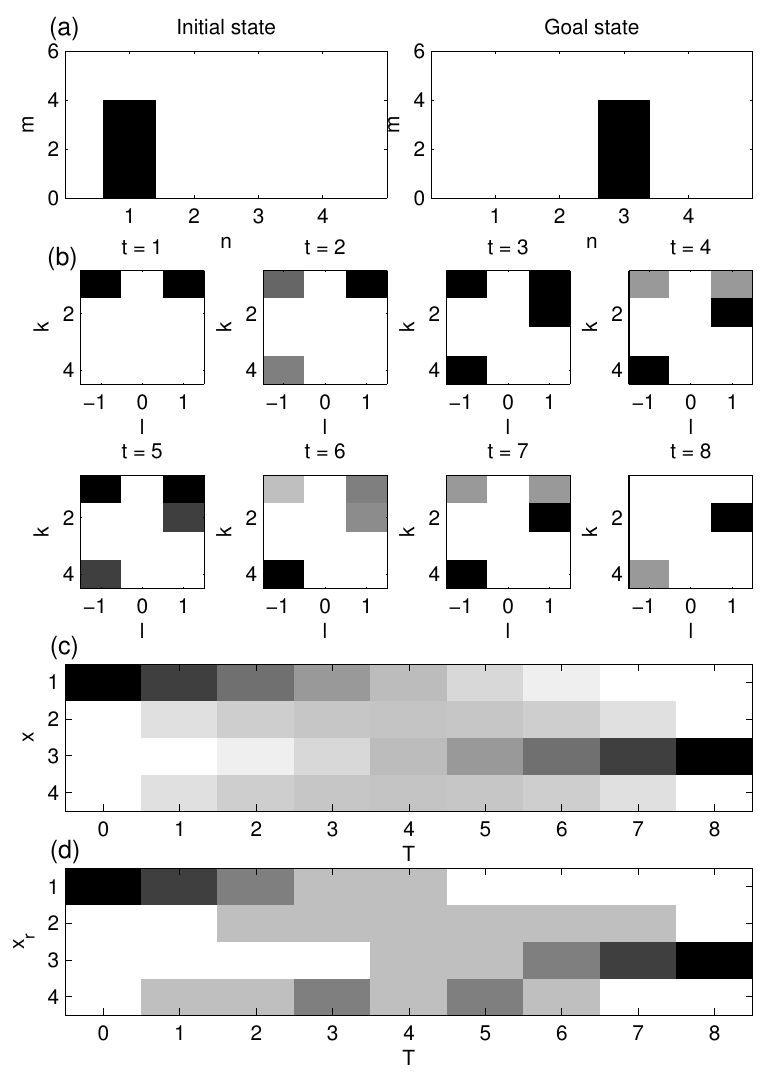}
\caption{
Control for the KL-blocks-world problem with end-cost:
example with $m=4, n=4$ and $T=8$.  \textbf{(a)} Initial and goal
states.  \textbf{(b)} Probability of action
$p(k^{t},l^{t}|x^{t-1})$ for each time step $t=1\ldots T$.
\textbf{(c)} Expected value $\langle x^t_i\rangle, i=1,\ldots,n$ given the
initial position and desired final position and \textbf{(d)}
the MAP solution for all times using a gray scale coding with white
coding for zero and darker colors coding for higher values.  
}
\label{fig:tower}
\ec
\end{figure}

Figure \ref{fig:tower} shows a small example where the planning task is to
shift a tower composed of four blocks which initially is at position $1$ to the
final position $3$.


To find the KL control we first condition the model both on the initial state
and the final state variables by ``clamping" all variables $x^1$ and $x^T$.
The KL control solution is obtained by computing for $t=1,\ldots,T$ the
marginal $p(k^t,l^t|x^{t-1})$.  In this case, we can find the exact solution
via the junction tree (JT) algorithm~\citep{lauritzen88,mooij2010libdai}.  The $k^t,l^t$ is
obtained by taking the MAP state of $p(k^t,l^t|x^{t-1})$ breaking ties at
random, which results in a new state $x_t$.

\begin{figure}
\bc
\includegraphics[width=0.6\textwidth]{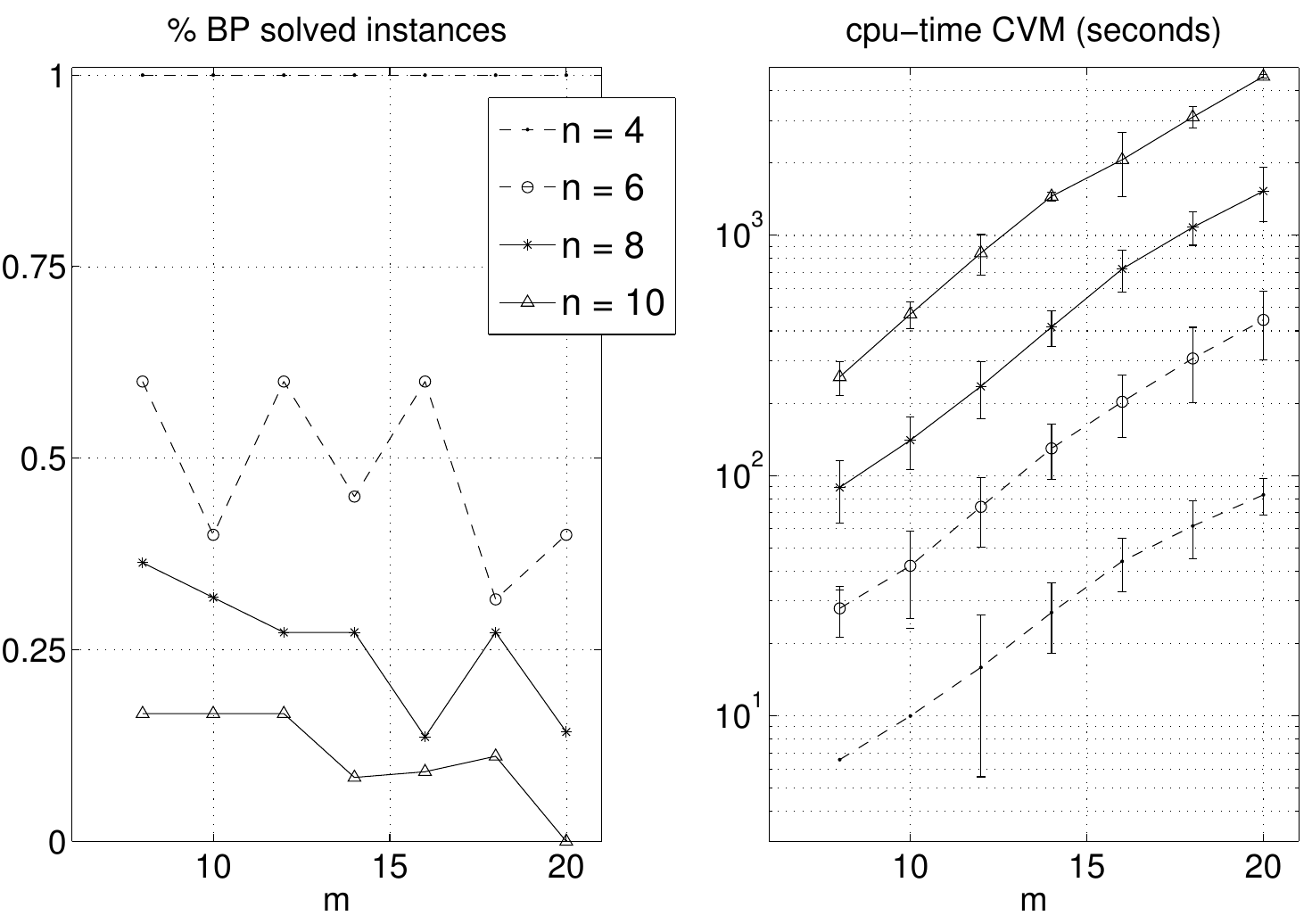}
\caption{
Control for the KL-blocks-world problem with end-cost:
results on approximate inference using random initial
and goal states.  \textbf{(Left)}~percent of instances where BP
converges for all $t=1:T$ as a function of $m$ for different values of
$n$.  \textbf{(Right)} CPU-time required for CVM to find
a correct plan for different values of $n,m$.  
$T$ was set to $\lceil \frac{m\cdot
n}{4} \rceil$.  We run $50$ instances for each pair $(m,n)$.
}\label{fig:stats-random}
\ec
\end{figure}

These probabilities $p(k^t,l^t|x^{t-1})$ are shown in Figure~\ref{fig:tower}b.  Notice that
the
symmetry in the problem is captured in the optimal control, which
assigns equal
probability when moving the first block to left or right
(Figure~\ref{fig:tower}b,c, t=1). Figure~\ref{fig:tower}d shows the
strategy resulting from the MAP estimate, which first unpacks the tower at position $1$ leaving all
four
locations with one block at $t=4$, and then re-builds it again at the goal
position $3$.


For larger instances, the JT method is not feasible because of too large tree
widths. For instance, to stack 4 blocks on 6 locations within a horizon of 11,
the junction tree has a maximal width of 12, requiring about 15 Gbytes of
memory.  We can nevertheless obtain approximate solutions using different
approximate inference methods.  In this work, we use the belief propagation
algorithm (BP) and a generalization known as the Cluster Variation method
(CVM).  We briefly summarize the main idea of the CVM method in
Appendix~\ref{appendix:cvm}.  We use the minimal cluster size, that is, the
outer clusters are equal to the interaction potentials $\psi$ as shown in the
graphical model Figure~\ref{graph}.

To compute the sequence of actions we follow again a sequential approach.
Figure \ref{fig:stats-random} shows results using BP and CVM. For $n=4$, BP
converges fast and finds a correct plan for all instances. For larger $n$, BP
fails to converge, more or less independently of $m$. Thus, BP can be applied
successfully to small instances only.  Conversely, CVM is able to find a
correct plan in all run instances, although at the cost of more CPU time, as
Figure \ref{fig:stats-random} shows.  The variance in the CPU error bars is
explained by the randomness in the number of actual moves required to solve
each instance, which is determined by the initial and goal states.
\begin{figure}[t]
\bc
\includegraphics[width=0.65\textwidth,angle=-90]{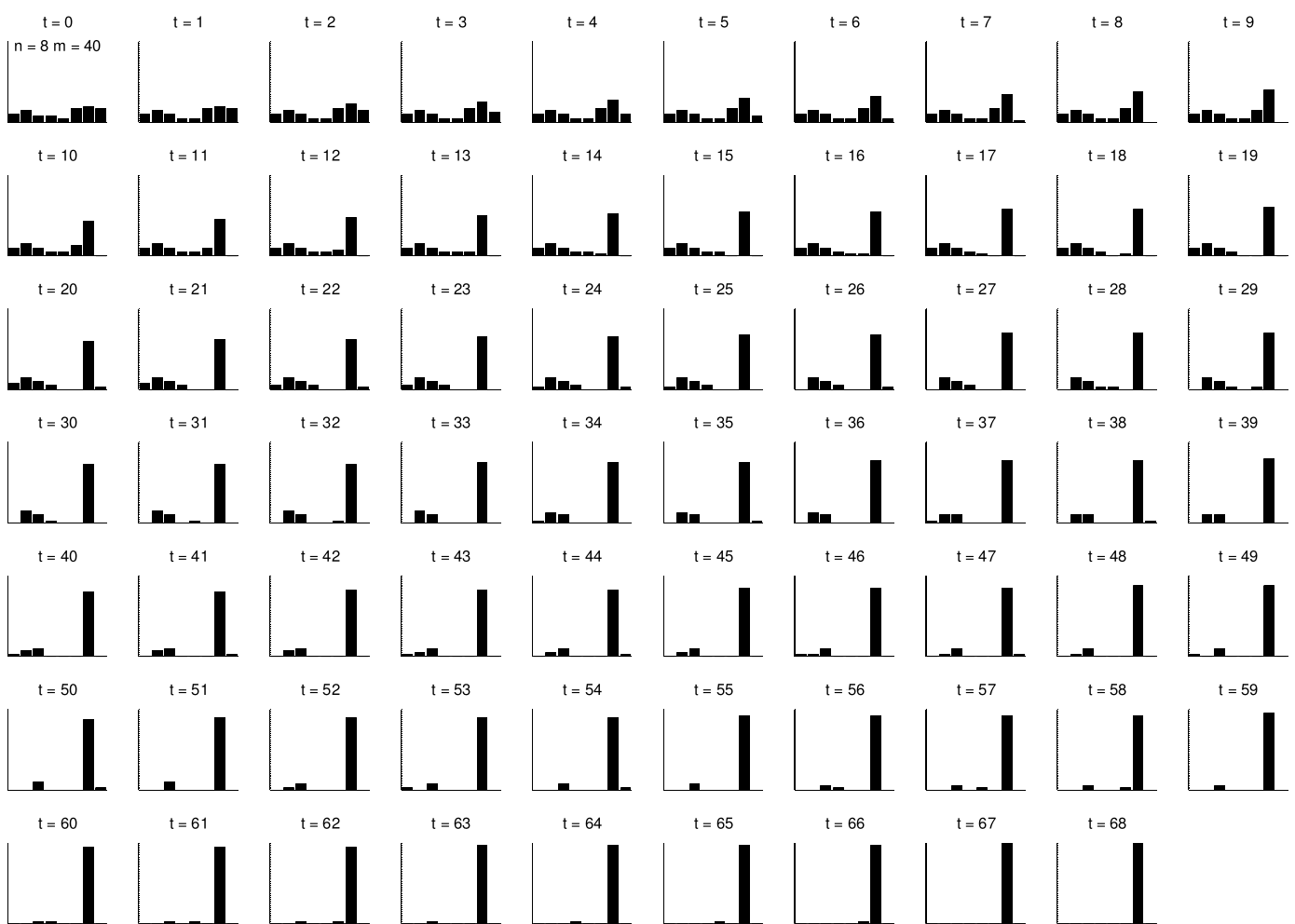}
\ec
\caption{Example of a large block stacking instance without end cost.  $n=8, m=40, T=80,
\lambda=10$ using CVM.} 
\label{large}
\end{figure}

\subsubsection{No goal state and $\lambda>0$: entropy minimization}
We now consider the problem without conditioning on $x^T$ and $\lambda>0$.
Although this may seem counter intuitive, removing the end constraint in fact
makes this problem harder, as the number of states that have significant
probability for large $t$ is much larger. BP is not able to produce any
reliable result for this problem.  We applied CVM to a large block stacking
problem with $n=8, m=40, T=80$ and $\lambda=10$.  We use again the minimal
cluster size and the double loop method of \cite{hes03a}.  The results are
shown in Figure \ref{large}.

The computation time was approximately 1 hour per $t$ iteration and memory use
was approximately 27 Mb. This instance was too large to obtain exact results.
We conclude that, although the CPU time is large, the CVM method is capable to
yield an apparently accurate control solution for this large
instance.

\section{Multi Agent cooperative game (KL-stag-hunt)}
In this section we consider a variant of the stag hunt game, a prototype game
of social conflict between personal risk and mutual benefit \citep{staghunt}.
The original two-player stag hunt game proceeds as follows: there are two
hunters and each of them can choose between hunting hare or hunting stag,
without knowing in advance the choice of the other hunter.  The hunters can
catch a hare on their own, giving them a small reward.  The stag has a much
larger reward, but it requires both hunters to cooperate in catching it.
 \begin{table}[h]
 \caption{Two-player stag hunt payoff matrix example: rows and columns indicate
 actions of one and the other player respectively. The payoff describes the
 reward for each hunter. For instance, if both go for the stag, they both get a
 reward of $3$. If one hunter goes for the stag and the other for the
 hare, they get a reward of $0$ and $1$ respectively.  }
 \label{tb:payoff}       
 \begin{tabular}{c|c|c|}
 \noalign{\smallskip}
 &Stag & Hare \\
 \hline
 Stag & $\mathbf{3,3}$ & $0,1$ \\
 Hare & $1,0$ & $\mathbf{1,1}$ \\
 \hline
 \end{tabular}
 \end{table}

Table \ref{tb:payoff} displays a possible payoff matrix for a stag hunt game.
It shows that both stag hunting and hare hunting are \emph{Nash equilibria},
that is, if the other player chooses stag, it is best to choose stag
(\emph{payoff} equilibrium, top-left), and if the other player chooses hare, it
is best to choose hare (\emph{risk-dominant} equilibrium, bottom-right).  It is
argued that these two possible outcomes makes the game socially more
interesting, than for example the \emph{prisoners dilemma}, which has only one Nash
equilibrium. The stag hunt allows for the study of cooperation within social
structures \citep{skyrms} and for studying the collaborative behavior of
multi-agent systems \citep{plos}.

We define the KL-stag-hunt game as a multi-agent version of the original stag
hunt game where $M$ agents live in a grid of $N$ locations and can move to
adjacent locations on the grid.  The grid also contains $H$ hares and $S$ stags
at certain fixed locations.  Two agents can cooperate and catch a stag
together with a high payoff $R_s$. Catching a stag with more than two agents
is also possible, but it does not increase the payoff. The agents can also
catch a hare individually, obtaining a lower payoff $R_h$.  The game is played
for a finite time $T$ and at each time-step all the agents perform an action.
The optimal strategy is thus to coordinate pairs of agents to go for different
stags.

Formally, let $x_i^t= 1,\hdots,N, i=1,\hdots,M, t=1,\hdots,T$ denote the position of
agent $i$ at time $t$ on the grid. Also, let $s_j=1,\hdots,N,
j=1,\hdots,S$, and $h_k= 1,\hdots,N, k=1,\hdots,H$ denote
the positions of the $j$th stag and the $k$th hare respectively.
We define the following state dependent reward as:
\begin{align*}
R(x^t) & =
    R_h\sum_{k=1}^H\sum_{i=1}^M \delta_{x_i^t,h_k}+
    R_s\sum_{j=1}^S\mathcal{I}\left\{\left( \sum_{i=1}^M x_i^t=s_j\right) > 1\right\},
\end{align*}
where $\mathcal{I}\{\cdot\}$ denotes the indicator function.  The first term
accounts for the agents located at the position of a hare.  The second one
accounts for the rewards of the stags, which require that at least two agents
to be on the same location of the stag.  Note that the reward for a stag is not
increased further if more than two agents go for the same stag.  Conversely, the
reward corresponding to a hare is proportional to the number of agents at its
position. 

The uncontrolled dynamics factorizes among the agents.  It allows an agent to
stay on the current position or move to an adjacent position (if possible) with
equal probability, thus performing a random walk on the grid.  Consider the
state variables of an agent in two subsequent time-steps expressed in Cartesian
coordinates, $x_i^t = \langle l,m\rangle, x_i^{t+1}=\langle l',m'\rangle$.  We
define the following function:
\begin{align*}
\psi_q\left( \langle l',m'\rangle, \langle l,m\rangle\right):=\mathcal{I} \Bigl\{&
    \left( (l'=l)   \wedge (m'=m) \right) \vee  \notag\\
&    \left( (l'=l-1) \wedge (m'=m) \wedge (l>0) \right) \vee \notag\\
&    \left( (l'=l) \wedge (m'=m-1) \wedge (m>0) \right) \vee \notag\\
&    \left( (l'=l+1) \wedge (m'=m) \wedge (l<\sqrt{N}) \right) \vee \notag\\
&    \left( (l'=l) \wedge (m'=m+1) \wedge (m<\sqrt{N}) \right)\Bigr\},
\end{align*}
that evaluates to one if the agent does not move (first condition), or if it
moves left, down, right, up (subsequent conditions) inside the grid boundaries.
The uncontrolled dynamics for one agent can be written as conditional
probabilities after proper normalization:
\begin{align*}
q\left( x_i^{t+1}= \langle l',m'\rangle| x_i^t=\langle l,m\rangle\right) &=
\frac
{\psi_q( \langle l',m'\rangle, \langle l,m\rangle )}
{\sum_{a,b}{\psi_q( \langle a,b\rangle, \langle l,m\rangle)  }}.
\end{align*}
and the joint uncontrolled dynamics become: 
\begin{align*}
q(x^{t+1}|x^{t}) & = \prod_{i=1}^M q(x_i^{t+1}|x_i^t).
\end{align*}
Since we are interested in the final configuration at end time $T$, we set the
state dependent path cost to zero for $t=1,\hdots,T-1$ and to
$\exp\left(-\frac{1}{\lambda}R(x^T)\right)$ for the end time.
\begin{figure}[h]
\bc
\includegraphics[angle=-90,width=.65\textwidth]{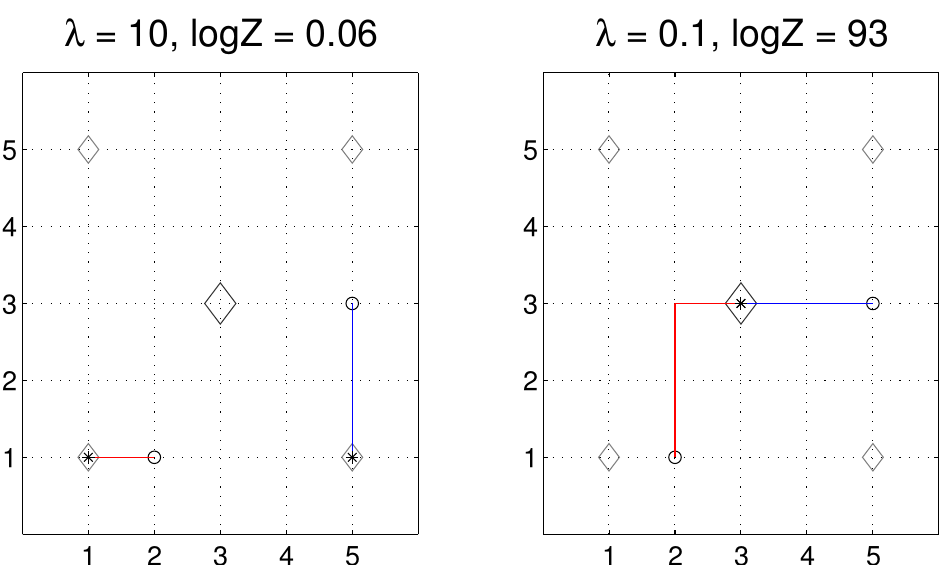}
\ec
\caption{Exact inference KL-stag-hunt:
Two hunters in a small grid.  There are four hares at each corner of the grid
(small diamonds) and one stag in the middle (big diamond).  Initial positions
of the hunters are denoted by small circles. One hunter is close to a hare and
the other is at the same distance of the stag and two hares. Final positions are
denoted by asterisks. The optimal paths are drawn in blue and red (Color online).
(\textbf{Left}) For $\lambda=10$, the optimal control is risk dominant, and
hunters go for the hares. (\textbf{Right}) For $\lambda=0.1$, the payoff
dominant control is optimal and hunters cooperate. $N=25, T=4, R_s = -10$ and
$R_h=-2$.
}
\label{fig:exact}
\end{figure}

To minimize $C$ in Equation \eqref{C}, exact inference in the joint
space can be done by backward message passing, using the following equations:
\begin{align}\label{eq:beta}
\beta^t({x}^t) & =\begin{cases}
 \exp\left(-\displaystyle\frac{1}{\lambda}R({x}^t)\right)                                       & \text{for $t=T$}\\
 \displaystyle\sum_{{x}^{t+1}} q({x}^{t+1}|{x}^{t})\beta({x}^{t+1}) & \text{for $t<T$}
 \end{cases},
 \end{align}
and the desired marginal probabilities can be obtained from the $\beta$-messages:
\begin{align}\label{eq:pbeta}
p({x}^{t+1}|{x}^{t}) & \propto q({x}^{t+1}|{x}^{t})\beta({x}^{t+1}).
\end{align}

To illustrate this game, we consider a small $5\times 5$ grid with two hunters
and apply equations \eqref{eq:beta} and \eqref{eq:pbeta}.  There are four hares
at each corner of the grid and one stag in the middle.  The initial positions
of the hunters are selected in a way that one hunter is close to a hare and the
other is at the same distance of the stag and two hares.  Starting from the
initial fixed state $x^1$, we select the next state according to the most
probable state from $p(x_i^{t+1}|x_i^t)$ until the end time.  We break ties
randomly.  Figure \ref{fig:exact} shows one resulting trajectory for two values
of $\lambda$.

For high values of $\lambda$ (left plot), each hunter catches one of the hares.
In this case, the cost function is dominated by KL term.  For small enough
values of $\lambda$ (right plot), both hunters cooperate to catch the stag. In
this case, the state cost, function $R(x^T)$, governs the optimal control cost.
Thus $\lambda$ can be seen as a parameter that controls whether the optimal
strategy is risk dominant or payoff dominant.

Note that computing the exact solution using this procedure becomes infeasible
even for small number of agents, since the joint state space scales as $N^M$.
In the next section, we show a more efficient representation using a factor
graph for which approximate inference is tractable.
\subsection{Graphical model for the KL-stag-hunt game}
\begin{figure}[t]
\bc
\includegraphics[width=.8\textwidth]{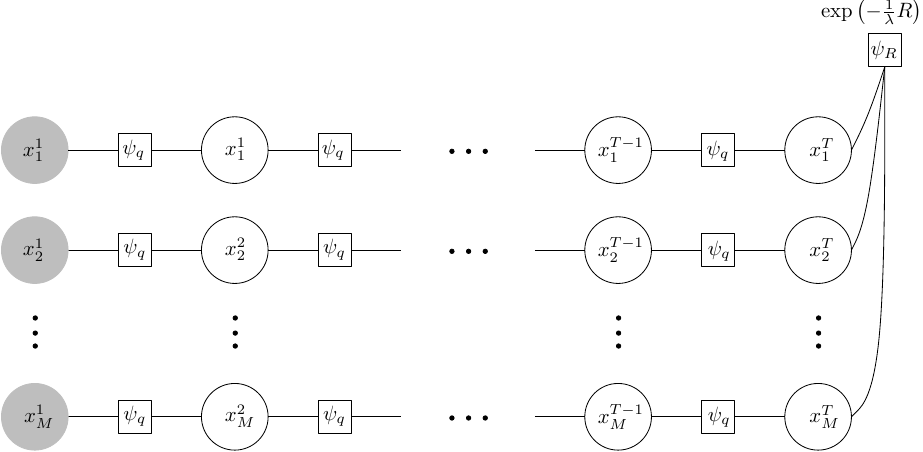}
\ec
\caption{Factor graph representation of the KL-stag-hunt problem.  Circles
denote variable nodes (states of the agents at a given time-step) and squares
denote factor nodes. There are two types of factor nodes: the ones
corresponding to the uncontrolled dynamics $\psi_q$ and the one corresponding
to the state cost $\psi_R$. Initial configuration in gray denotes the states
``clamped" to an initial given value. Despite being a tree, exact inference and
approximate inference are intractable in this model due to the complex factor
$\psi_R$.
}
\label{fig:MASfg}
\end{figure}

The corresponding graphical model of the KL-stag-hunt game is depicted in
Figure~\ref{fig:MASfg} as a factor graph.  Since the uncontrolled dynamics
factorizes over the agents, the joint state can be split in different 
variable nodes. 
Note that since there is only state cost at the end time, the graphical model
becomes a tree. However, the factor node associated to the state cost function
$\psi_R(x^T):=\exp(-\frac{1}{\lambda} R(x^T))$ involves all the agent states,
which still makes the problem intractable. Even approximate inference
algorithms such as BP can not be applied, since messages from $\psi_R$ to one
of the state variables $x_i^T$ would require a marginalization involving a sum
of $(N-1)^M$ terms.

\begin{figure}
\bc
\includegraphics[width=.55\textwidth]{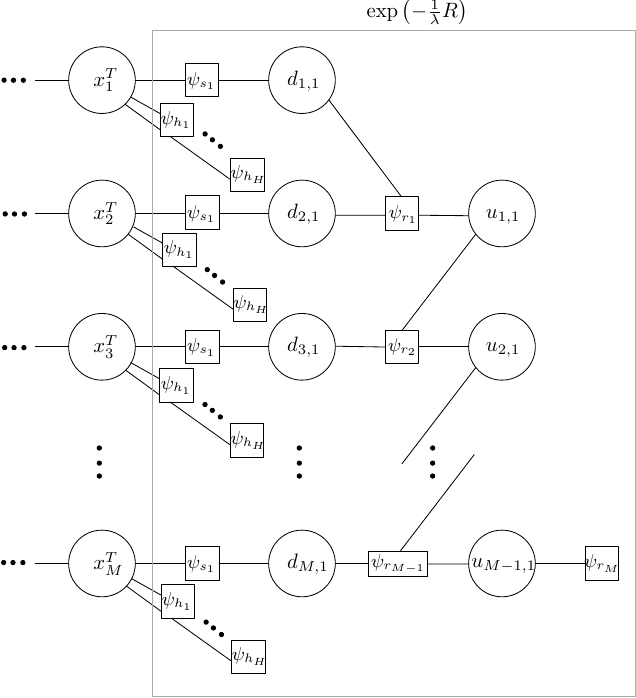}
\ec
\caption{Decomposition of the complex factor
$\psi_R$
into simple factors involving at most three variables of small cardinality.
Each state variable is linked to $H$ factors corresponding to the hares
locations.  For each stag there is a chain of factors $\psi_{r_i}, 
i=1,\hdots,M-1$ which evaluates to one for the allowed configurations and to
zero otherwise. Factor $\psi_{r_M}$ weights the configuration of having zero,
one or more agents being at the stag position (figure shows the case of one
stag only).  }
\label{fig:MAStrick}
\end{figure}

However, we can exploit the particular structure of that factor by decomposing
it in smaller factors defined on small sets of (at most three) auxiliary variables
of small cardinality.  This transformation becomes intuitive once the graphical
model representation for the problem is identified.  The procedure defines
indicator functions for the allowed configurations which are weighted according
to the corresponding cost. Figure \ref{fig:MAStrick} illustrates the procedure
for the case of one stag.
\begin{enumerate}
\item First, we add $H\times M$ factors $\psi_{h_k}(x_i^T)$, defined for each
hare location $h_k$ and each agent variable $x_i^T$. These factors account for the hare costs:
      \begin{align*}
      \psi_{h_k}(x_i^T) := &
      \begin{cases}
      \exp\left(-\frac{1}{\lambda}R_h\right) & \text{if $(x_i^T = h_k)$ }\\
      1                                    & \text{otherwise}.
      \end{cases}
      \end{align*}
\item Second, we add factors $\psi_{s_j}(x_i^T,d_{i,j})$ for
each stag $j$ defined on
each state variable $x_i^T$ and new introduced binary variables $d_{i,j}=0,1$.
These factors evaluate to one when variable $d_{i,j}$ takes the value
of a Kronecker $\delta$ of the agent's state $x_i^T$ and the position
of a stag $s_j$, and zero otherwise:
\begin{align*}
\psi_{s_j}(x_i^T,d_{i,j}) := \mathcal{I}\Bigl\{(d_{i,j} & = \delta_{x_i^T,s_j})\Bigr\}.
\end{align*}
\item
Third, for each stag, we introduce a chain of factors that involve the binary
variables $d_{i,j}$ and additional variables $u_{i,j}=0, 1, 2$.  The new
variables $u_{i,j}$ encode whether the stag $j$ has zero, one, or more agents
after considering the $(i+1)$th agent.  The new factors are: 
\begin{align*}
\psi_{r_1}(d_{1,j}, d_{2,j}, u_{1,j}) := 
    \mathcal{I}\Bigl\{ & \left( (d_{1,j}=0) \wedge (d_{2,j}=0) \wedge (u_{1,j}=0) \right) \vee \\
    & \left( (d_{1,j}=1) \wedge (d_{2,j}=1) \wedge (u_{1,j}=2) \right) \vee\\
    & \left( (d_{1,j} \neq d_{2,j})   \wedge (u_{1,j}=1) \right)
    \Bigr\}. \\
\psi_{r_{i-1}}(u_{i-1,j}, d_{i,j}, u_{i,j}) := 
    \mathcal{I}\Bigl\{ &
      \left( (d_{i,j}=0) \wedge (u_{i-1,j} = u_{i,j}) \right) \vee \\
    & \left( (d_{i,j}=1) \wedge (u_{i-1,j} = 0) \wedge (u_{i,j}=1) \right)\vee \\
    & \left( (d_{i,j}=1) \wedge (u_{i-1,j} = 1) \wedge (u_{i,j}=2) \right)\vee \\
    & \left( (d_{i,j}=1) \wedge (u_{i-1,j} = 2) \wedge (u_{i,j}=2) \right)
    \Bigr\}. \\
\end{align*}
\item Finally, we define factors $\psi_{r_M}$ that weight the allowed configurations:
\begin{align*}
\psi_{r_M}(u_{M-1,j}) := & 
    \begin{cases}
       \exp\left(-\frac{1}{\lambda}R_s\right) & \text{if $(u_{M-1,j} = 2)$ }\\
       1                                      & \text{otherwise}.
    \end{cases}
\end{align*}
\end{enumerate}
The original factor can be rewritten marginalizing the 
auxiliary variables $d_{i,j}, u_{i,j}$ over the product of the previous factors
$\psi_{s_j},\psi_{h_k},\psi_{r_i}$:
\begin{align*}
\exp\left(-\frac{1}{\lambda}R(x^T)\right) & = \psi_S(x^T)\psi_H(x^T)\\
\psi_S(x^T) & := \prod_{j=1}^S\Bigl[
\sum_{
\substack{d_{1,j},d_{2,j}\\
u_{1,j},u_{M-1,j}}}
\left( \psi_{s_j}(x_1^T,d_{1,j})\psi_{s_j}(x_2^T,d_{2,j} )  \right)
\psi_{r_1}(d_{1,j},d_{2,j},u_{1,j})\\
&\psi_{r_M}(u_{M-1,j})
\sum_{\substack{
d_{3,j},\hdots,d_{M,j}\\
u_{2,j},\hdots,u_{M,j}}}
\prod_{i=3}^M\psi_{r_{i-1}}(u_{i-1,j},d_{i,j},u_{i,j})\psi_{s_j}(x_i^T,d_{i,j})\Bigr]\\
\psi_H(x^T) & := \prod_{k=1}^H \psi_{h_k}(x_i^T),
\end{align*}
where for clarity of notation we have grouped the factors related to the stags
and hares in $\psi_S(x^T)$ and $\psi_H(x^T)$, respectively.
\begin{figure}[!t]
\bc
\includegraphics[angle=-90,width=.75\textwidth]{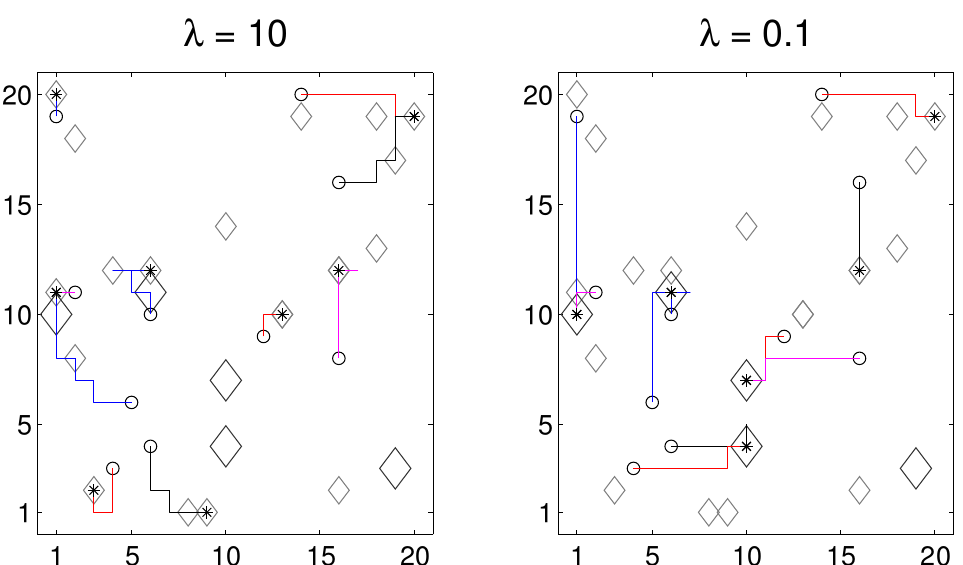}
\ec
\caption{Approximate inference KL-stag-hunt: Control obtained using BP for
$M=10$ hunters in a large grid.  See Figure \ref{fig:exact} for a description
of the symbols. (\textbf{Left}) Risk dominant control is obtained for
$\lambda=10$, where all hunters go for a hare. (\textbf{Right}) Payoff dominant
control is obtained for $\lambda=0.1$.  In this case, all hunters cooperate to
capture the stags except the ones on the upper-right corner, who are too far
away from the stag to reach it in $T=10$ steps. Their optimal choice is to go
for a hare. $N=400, S = M/2, R_s = -10, H=2M$ and $R_h=-2$.
}
\label{fig:klstaghunt_approx}
\end{figure}

The extended factor graph is tractable since it involves factors of no
more than three variables of small cardinality. Note that this transformation
can also be applied if additional state costs are incorporated at each time-step
$\psi_R(x^t)\neq 0, t = 1,\hdots,T$. However, such a representation is not of
practical interest, since it complicates the model unnecessarily.

Finally, note that the tree-width of the extended graph still grows fast with
the number of agents $M$ because variables $d_{i,j}$ and $u_{i,j}$ are coupled.
Thus, exact inference using the JT algorithm is still possible on small
instances only.
\subsection{Approximate inference of the KL-stag-hunt problem}
In this section we analyze large systems for which exact inference is not
possible using the JT algorithm.  The belief propagation (BP) algorithm is an
alternative approximate algorithm that we can run on the previously described
extended factor graph. 

We use the following setup: for a fixed number of agents $M$, we set the number
of stags $H=2M$ and the number of hares $S=\frac{M}{2}$.  Their locations, as
well as the initial states $x^1$ are chosen randomly and non-overlapping.
We then construct a factor graph with initial states
``clamped" to $x^1$ and build instance-dependent factors $\psi_{s_j}$ and
$\psi_{h_k}$.  We run BP using sequential updates of the messages.  If BP
converges in less than $500$ iterations, the optimal trajectories of the agents
are computed using the estimated marginals (factor beliefs) for
$\psi_q(x^{t+1}|x^t)$ after convergence. Starting from $x^1$, we select the
next state according to the most probable state from $p_{BP}(x_i^{t+1}|x_i^t)$
until the end time.  We break ties randomly.
We analyze the system as a function of parameter $\lambda$ for a several number
of realizations.

The global observed behavior is qualitatively similar to the
one of a small system: for $\lambda$ very large, a risk-dominant control is
obtained and for $\lambda$ small enough, payoff control dominates. This is
behavior is illustrated in Figure \ref{fig:klstaghunt_approx}, where an
example for $\lambda=10$ and $\lambda=0.1$ are shown.  We can thus conclude
that BP provides an efficient and good approximation for large systems where
exact inference is not feasible. 
\begin{figure}[!t]
\bc
\includegraphics[angle=-90,width=.7\textwidth]{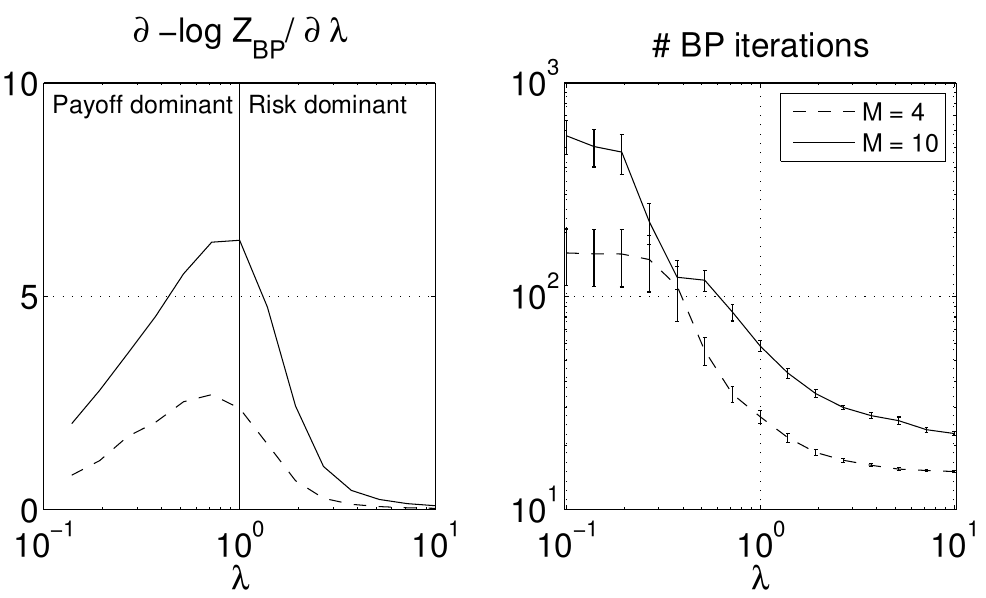}
\ec
\caption{Approximate inference KL-stag-hunt:
(\textbf{Left})
Change in the expected cost with respect to $\lambda$ as a function of
$\lambda$ for $\langle M = 4, N=100\rangle$ and $\langle M = 10, N=225\rangle$.
The curve becomes sharper and its maximum gets closer to $\lambda=1$ for larger
systems, suggesting a phase transition phenomenon between the risk dominant and
the payoff dominant regimes.  (\textbf{Right}) Number of BP iterations required
for convergence as a function of $\lambda$.  Results are averages over $20$
runs with random initial states.  $R_s = -10, R_h=-2$ and $T=10$.  }
\label{fig:bp_klstaghunt}
\end{figure}

To characterize the solutions, we compute the approximated expected cost as in
\eqref{path-integral}, that is $-\log Z_{BP}$.  We observe that for large
systems that quantity changes abruptly at $\lambda\approx 1$.  Qualitatively,
the optimal control obtained on the boundary between risk-dominant and
payoff-dominant strategies differs maximally between individual instances and
strongly depends on the initial configuration.  This suggests a phase
transition phenomenon typical of complex physical systems, in this case
separating the two kind of optimal behaviors, where $\lambda$ plays the role of
a ``temperature" parameter.

Figure \ref{fig:bp_klstaghunt} shows this effect. The left plot shows the
derivative of the expected approximated cost averaged over $20$ instances.  The
curve becomes sharper and its maximum gets closer to $\lambda=1$ for larger
systems. Error bars of the number of iterations required for convergence is
shown on the right.  The number of BP iterations quickly increases as we
decrease $\lambda$, indicating that the solution for which agents cooperate is
more complex to obtain. For $\lambda$ very small, BP may fail to converge
after $500$ iterations.

\section{Related work}
\label{related_work}

The idea to treat a control problem as an inference problem has a long history.
The best known example is the linear quadratic control problem, which is
mathematically equivalent to an inference problem and can be solved as a Kalman
smoothing problem \citep{stengel93}.  The key insight is that the value
function that is iterated in the Bellman equation becomes the (log of the)
backward message in the Kalman filter.  The exponential relation was
generalized in \citet{kappen_prl05} for the non-linear continuous space and
time (Gaussian case) and in \citet{NIPS2006_691} for a class of discrete
problems.

There is a line of research on how to compute optimal action sequences in
influence diagrams using the idea of probabilistic inference
\citep{cooper88,tatman_shachter1990,shachter_peot1992}. Although this technique
can be implemented efficiently using the junction tree approach for single
decisions, the approach does not generalize in an efficient way to optimal
decisions, in the expected-reward sense, in multi-step tasks.  The reason is
that the order in which one marginalizes and optimizes strongly affects the
efficiency of the computation.  For a Markov decision process (MDP) there is an
efficient solution in terms of the Bellman equation \footnote{Here we mean by
efficient, that the sum or min over a sequence of states or actions can be
performed as a sequence of sums or mins over states.}. For a general influence
diagram, the marginalization approach as proposed in
\citet{cooper88,tatman_shachter1990,shachter_peot1992} will result in an
intractable optimization problem over $u^{0:T-1}$ that cannot be solved
efficiently (using dynamic programming), unless the influence diagram has an
MDP structure.

The KL control theory shares similarities with work in reinforcement learning
for policy updating.  The notion of KL divergence appears naturally in the work
of \citet{bagnell_schneider2003} who proposes an information geometric approach
to compute the natural policy gradient (for small step sizes).  This idea is
further developed into an Expectation-Maximization (EM) type algorithm
\citep{dayan_hinton1997} in recent work
\citep{peters_muelling_altun2010,kober_peters2011} using a relative entropy
term. The KL divergence acts here as a regularization that weights the relative
dependence of the new policy on the data observed and the old policy,
respectively.

It is interesting to compare the the notion of free energy in continuous-time
dynamical systems with Gaussian noise considered in \citet{Friston} with the
path integral formalism of \citet{kappen_prl05}, which is a special case of KL
control theory.  \citet{Friston} advocate the optimization of free energy as a
guiding principle to describe behavior of agents. The main difference between
the KL control theory and Friston's free energy principle is that in KL control
theory, the KL divergence plays the role of an expected future cost and its
optimization yields a (time dependent) optimal control trajectory, whereas
Friston's free energy computes a control that yields a time-independent
equilibrium distribution, corresponding to the minimal free energy. Friston's
free energy formulation is obtained as a special case of KL control theory when
the dynamics and the reward/cost is time-independent and the horizon time is
infinite.

The KL control approach proposed in this paper also bears some relation to the
EM approach of \cite{toussaint_icml2006}, who consider the discounted reward
case with 0,1 rewards.  The posterior can be considered a mixture over times at
which rewards are incorporated.  For an homogeneous Markov process and time
independent costs, the backward message passing can be effectively done in a
single chain and not the full mixture distribution needs to be considered.  We
can compare the EM approach of \cite{toussaint_icml2006} (TS) and the KL
control approach (KL):
\begin{itemize}
\item The TS approach is more general than the KL approach, in the sense that
the reward considered in TS is an arbitrary function of state and action
$R(x,u)$, whereas the reward considered in KL is a sum of a state dependent
term $R(x)$ and a KL divergence.
\item 
The KL approach is significantly more efficient than the TS approach.
In the TS approach, the backward messages are computed for a fixed
policy $\pi$ (E-step), from which an improved policy is computed
(M-step). This procedure is iterated until convergence. In the KL
approach, the backward messages give the optimal control directly,
with no further need for iteration. 
\item In addition, the KL approach is more efficient than the TS approach for
time-dependent problems.
Using the TS approach for time-dependent problems makes the computation
a factor $T$ more time-consuming than for the time-independent case,
since all mixture components must be computed. The complexity of
the KL control approach does not depend on whether the problem is
time-dependent or not.
\item The TS and KL approach optimize with respect to a different
quantity.  The TS approach writes the state transition $p(y|x) =
\sum_u p(y|x,u) \pi(u|x)$ and optimizes with respect to $\pi$. The
KL approach optimizes the state transition probability $p(y|x)$
directly either as a table or in a parametrized way.
\end{itemize}

\section{Discussion}
In this paper, we have shown the equivalence of a class of stochastic optimal
control problems to a graphical model inference problem. As a result, exact or
approximate inference methods can directly be applied to the intractable stochastic
control computation.  The class of KL control problems contains interesting
special cases such as the continuous non-linear Gaussian stochastic control
problems introduced in \citet{kappen_prl05}, discrete planning tasks and
multi-agent games, as illustrated in this paper.

We notice, that there exist many stochastic control problems that are
outside of this class.  In the basic formulation of Equation \eqref{Cmdp},
one can construct control problems where the functional form of the
controlled dynamics $p^t(x^{t+1}|x^t,u^t)$ is given as well as the cost
of control $R(x^t,u^t,x^{t+1},t)$. In general, there may then not exist
a $q^t(x^{t+1}|x^t)$ such that Equation \eqref{Rspecial} holds.

In this paper, we have considered the model based case only. The extension to
the model free case would require a sampling based procedure. See \cite{joris}
for initial work in this direction.

We have demonstrated the effectiveness of approximate inference methods to
compute the approximate control in a block stacking task and
a multi-agent cooperative task.

For the KL-blocks-world, we have shown that an entropy minimization task is
more challenging than stacking blocks at a fixed location (goal state), because
the control computation needs to find out where the optimal location is.
Standard BP does not give any useful results if no goal state was
specified, but apparently good optimal control solutions were obtained
using generalized belief propagation (CVM).  We found that the marginal
computation using CVM is quite difficult compared to other problems that have
been studied in the past \citep{albers_genetics07}, in the sense that relatively
many inner loop iterations were required for convergence.  One can  improve the
CVM accuracy, if needed, by considering larger clusters
\citep{Yedidia05constructingfree} as has been demonstrated in other contexts
\citep{albers_bmc05}, at the cost of more computational complexity.  

We have given evidence that the KL control formulation is particularly
attractive for multi-agent problems, where $q$ naturally factorizes over agents
and where interaction results from the fact that the reward depends on the
state of more than one agent.  A first step in this direction was already made
in \citet{wimw_uai06,broek2006}.  In this case, we have considered the
KL-stag-hunt game and shown that BP provides a good approximation and allows to
analyze the behavior of large systems, where exact inference is not feasible. 
 
We found that, if the game setting strongly penalizes large deviations from the
baseline (random) policy, the coordinated solution is sub-optimal. That means
that the optimal solution distributes the agents among the different hares
rather than bringing them jointly to the stags (risk-dominant regime).  On the
contrary, if the agents are not constrained by deviating too much from the
baseline policy to maximize $\langle R\rangle$, the coordinated solution
becomes optimal (payoff dominant regime). We believe that this is an
interesting result, since it provides a explanation of the emergence of
cooperation in terms of an effective temperature parameter $\lambda$.


\subsubsection*{Acknowledgments}
We would like to thank anonymous reviewers for helping on improving the
manuscript, Kees Albers for making available his sparse CVM code, Joris Mooij
for making available the libDAI software and Stijn Tonk for useful
discussions.  The work was supported in part by the ICIS/BSIK consortium.


\appendix

\section{Boltzmann distribution}
\label{appendix:bm}
Consider the KL divergence between a normalized probability distribution $p(x)$ and some
positive function $\psi(x)$:
\beaa
C(p)=\sum_x p(x)\log \frac{p(x)}{\psi(x)}
\eeaa
$C$ is a function of the distribution $p$. We compute the distribution that minimizes $C$ with
respect to $p$ subject to normalization $\sum_x p(x)=1$ by adding a Lagrange multiplier:
\beaa
L(p)&=&C(p)+\beta \left(\sum_x p(x)-1\right)\\
\frac{d L}{d p(x)}&=& \log \frac{p(x)}{\psi(x)}+1+\beta
\eeaa
Setting the derivative equal to zero yields $p(x)=\psi(x)\exp(-\beta -1)=\psi(x)/Z$, where we
have defined $Z=\exp(\beta+1)$. The normalization
condition $\sum_x p(x)=1$ fixes $Z=\sum_x \psi(x)$. Substituting the solution for $p$ in the
cost $C$ yields $C=-\log Z$.

\section{Relation to continuous path integral model}
\label{appendix:continuous}
We write $p(x'|x)=\mathcal{N}(x'|x+f(x,t)dt +g(x,t)u(x,t) dt ,\Xi dt )$
with $\Xi(x,t)= g(x,t) \nu g(x,t)^T $ in Equation \eqref{controlleddynamics} as
\begin{align*}
p(x'|x)&=\cN(x'|x+f(x,t) dt, \Xi(x,t)dt )\exp\Bigl((\dot{x}-f(x,t))^T
\Xi^{-1}g(x,t) u(x,t)-\\
& \qquad\quad \frac{dt}{2}(g(x,t) u(x,t))^T \Xi^{-1}g(x,t) u(x,t) \Bigr)\\
&=q(x'|x)\exp\left(U(x,x',t) dt\right)\\
U(x,x',t)&=(\dot{x}-f(x,t))^T
\Xi^{-1}g(x,t)u(x,t)-\frac{1}{2}(g(x,t)u(x,t))^T \Xi^{-1} g(x,t)u(x,t) .
\end{align*}
with $\dot{x}=(x'-x)/dt$.

In order to make the link to Equation \eqref{C}
we compute
\begin{align*}
\sum_{x'}p(x'|x)\log \frac{p(x'|x)}{q(x'|x)}& = \sum_{x'}p(x'|x) U(x,x',t) dt\\
&=\frac{dt}{2}(g(x,t)u(x,t))^T \Xi(x,t)^{-1} g(x,t)u(x,t)\\
&=\frac{dt}{2}u(x,t)^T \nu^{-1} u(x,t),
\end{align*}
where we have made use of the fact that $\sum_{x'}p(x'|x) x'=x+f(x,t)dt +g(x,t)u(x,t) dt$ and 
$g^T \Xi^{-1} g= g^T (g^{-1})^T\nu^{-1} g^{-1} g=\nu^{-1}$. 
\footnote{
When $g$ is not a square matrix (when the number of controls is
less than the dimension of $x$), $g^{-1}$ denotes the pseudo-inverse of $g$.
For any $u$, the pseudo-inverse has the property that $g^{-1}g u=u$.}
Therefore,
\beaa
KL(p||q)&=&\sum_{x^{dt:T}}p(x^{dt:T}|x^0)\log\frac{p(x^{dt:T}|x^0)}{q(x^{dt:T}|x^0)}\\
&=&\sum_{s=0}^{T-dt}\sum_{x^s}p(x^s|x^0)\sum_{x^{s+dt}}p(x^{s+dt}|x^s)U(x^s,x^{s+dt},s)dt\\
&=&\sum_{s=0}^{T-dt}dt \sum_{x^s}p(x^s|x^0) \frac{1}{2}(u(x^s,s))^T \nu^{-1} u(x^s,s).
\eeaa

In the limit of $dt\rightarrow 0$ the KL divergence between $p$ and $q$ becomes 
\beaa
KL(p||q)&=&\av{\int_0^T dt \frac{1}{2}u(x(s),s)^T \nu^{-1} u(x(s),s)}
\eeaa
in agreement with Equation \eqref{cost_pi}.

\section{Cluster Variation Method}
\label{appendix:cvm}

In this appendix, we give a brief summary of the CVM method and the
double loop approach. For a more complete description see
\cite{yedidia00,kap01a,hes03a}.

\begin{figure}
\bc
\includegraphics[width=0.2\textwidth]{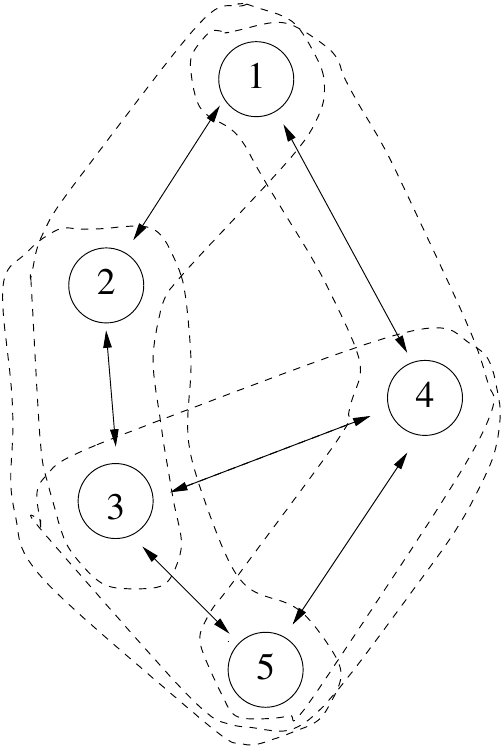}
\includegraphics[width=0.4\textwidth]{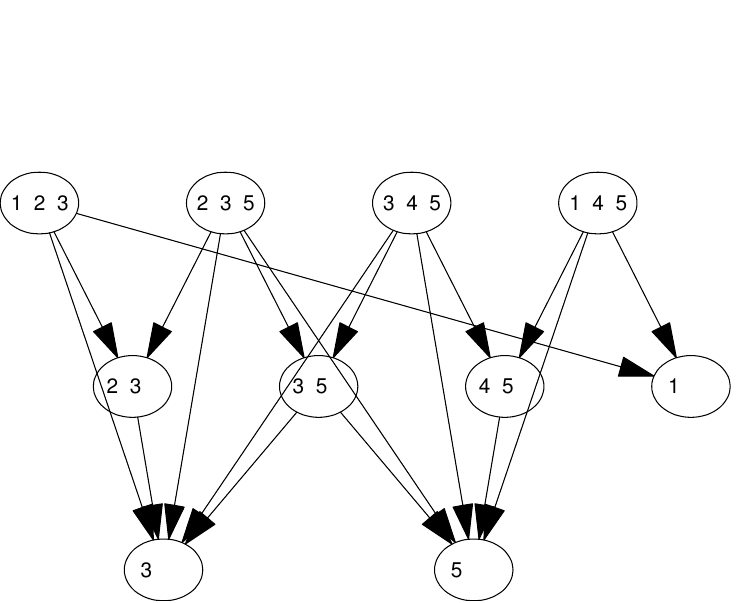}
\includegraphics[width=0.3\textwidth]{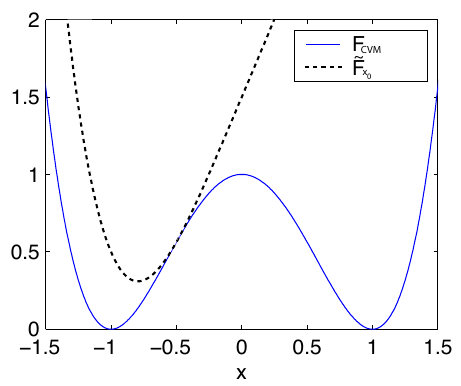}
\ec
\caption{\textbf{(Left)} Example of a small network and a choice of clusters
for CVM. \textbf{(Middle)} Intersections of clusters recursively define a set of
sub-clusters. \textbf{(Right)} $F_\mathrm{cvm}$ is non-convex (blue curve) and is
bounded by a
convex function $\tilde{F}_{x_0}$ (Color online). }
\label{cvm1}
\end{figure}

The cluster variation method replaces the probability distribution $p(x)$ in
the minimization Equation \eqref{C} by a large number of (overlapping)
probability distributions (clusters), each describing the interaction between a
small number of variables.
\begin{align*}
p(x) &\approx \{p_\alpha(x_\alpha),\alpha=1,\ldots\}
\end{align*}
with each $\alpha$ a subset of the indices $1,\ldots,n$, $x_\alpha$ the
corresponding subset of variables and $p_\alpha$ the probability distribution
on $x_\alpha$. The set of clusters is denoted by $B$, and must be such that any
interaction term $\psi_\alpha(x_\alpha)$, with $\psi(x)=\prod_\alpha
\psi_\alpha(x_\alpha)$ from Equation \eqref{psi}, is contained in at least one
cluster.  One denotes the set of all pairwise intersections of clusters in $B$,
as well as intersections of intersections by $M$.  Figure \ref{cvm1} (left) gives
an example of a small directed graphical model, where $B$ consists of 4
clusters and $M$ consists of 5 sub-clusters, Figure \ref{cvm1} (middle).

The CVM approximates the KL divergence, Equation \eqref{C}, as
\begin{align*}
C(x^0,p)&\approx F_\mathsf{cvm}(\{p_\alpha\})\\
F_\mathsf{cvm}(\{p_\alpha\})&=
\sum_{\alpha\in B} \sum_{x_\alpha}
p_\alpha(x_\alpha)\log\frac{p_\alpha(x_\alpha)}{\psi_\alpha(x_\alpha)}
+
\sum_{\beta \in M} a_\beta \sum_{x_\beta}
p_\beta(x_\beta)\log p_\beta(x_\beta).
\end{align*}

$F_\mathsf{cvm}$ is minimized with respect to all $\{p_\alpha\}$ subject to
normalization and consistency constraints:
\begin{align*}
\sum_{x_\alpha} p_\alpha(x_\alpha)=1,\qquad
p_\alpha(x_\beta)=p_\beta(x_\beta),\quad \beta\subset \alpha,\qquad
p_\alpha(x_\alpha)\ge 0
\end{align*}
The numbers $a_\beta$ are called the M\"obius or overcounting numbers.  They
can be recursively computed from the formula
\begin{align*}
1=\sum_{\alpha \in B\cup M, \alpha \supset \beta} a_\alpha, \qquad \forall
\beta \in B \cup M
\end{align*}

Since $a_\alpha$ can be both positive and negative, $F_\mathrm{cvm}$ is
not convex. A guaranteed convergent approach to minimize $F_\mathrm{cvm}$
is a double loop approach where the outer loop is to upper-bound
$F_\mathrm{cvm}$ by a convex function $\tilde{F}_{p^0}$ that
touches at the current cluster solution $p^0=\{p^0_\alpha\}$.  Optimizing $\tilde{F}_{p^0}(p)$
 is a convex problem that can be solved using the dual approach (inner
loop) and is guaranteed to decrease $F_\mathrm{cvm}$ to a local
minimum. 
The solution $p^*(p^0)$ 
of this convex sub-problem is guaranteed to decrease $F_\mathrm{cvm}$:
\[
F_\mathrm{cvm}(p^0)=\tilde{F}_{p^0}(p^0)\ge \tilde{F}_{p^0}(p^*(p^0)) \ge
F_\mathrm{cvm}(p^*(p^0))
\]
Based on $p^*(p_0)$ a new convex upper bound is computed (outer loop).
This is called a double loop method.
The approach is illustrated in Figure \ref{cvm1} (right).

Alternatively, one can choose to ignore the non-convexity issue.
Adding Lagrange multipliers $\lambda$ to enforce the constraints
one can minimize with respect to $p=\{p_\alpha\}$ and obtain an
explicit solution of $p$ in terms of the interactions $\psi$ and
the $\lambda$'s. Inserting this solution in the above constraints
results in a set of non-linear equations for the $\lambda$'s, which
one may attempt to solve by fixed point iteration. It can be shown
that these equations are equivalent to the message passing equations
of belief propagation. Unlike the above double loop approach, belief
propagation does not converge in general, but tends to give a fast and
accurate solution for those problems for which it does converge.

\bibliographystyle{apalike}

\bibliography{klcontrol}   


\end{document}